\documentclass[a4paper,12pt]{amsart}
\usepackage{latexsym,amsmath,amsthm,amssymb,amsfonts}
\usepackage[dvipsnames]{xcolor}
\usepackage{tikz}
\usepackage{geometry}
\numberwithin{equation}{section}
\usepackage{geometry}

\allowdisplaybreaks

\def\arcsinh{\operatorname{arcsinh}}

\def\s{\,\,\,\,}
\def\R{\mathbb{R}}

\def\C{\mathbb{C}}

\def\K{\mathbb{K}}

\def\endproof{$\hfill\Box$\\}

\def\vol{{\rm Area}}
\def\genus{\mathfrak{g}}

\newcounter{Cnumber}

\numberwithin{equation}{section}
\newtheorem{theorem}{Theorem}[section]
\newtheorem{lem}[theorem]{Lemma}
\newtheorem{thm}[theorem]{Theorem}
\newtheorem{pro}[theorem]{Proposition}
\newtheorem{cor}[theorem]{Corollary}
\newtheorem{defi}[theorem]{Definition}
\newtheorem{rem}[theorem]{Remark}

\makeatletter

\newcommand{\Rmnum}[1]{\expandafter\@slowromancap\romannumeral #1@}
\makeatother
\title[]{\bf 
Integral gradient estimates on a closed surface}

\author{Yuxiang Li, Rongze Sun}
\address{Department of Mathematical Sciences, Tsinghua University, People's Republic of China}
\email{liyuxiang@tsinghua.edu.cn,\s srz20@mails.tsinghua.edu.cn}
\date{}

\begin{document}

\begin{abstract}
Let $(\Sigma, g)$ be a closed Riemann surface, and let $u$ be a weak solution to equation
\[
- \Delta_g u = \mu,
\]
where $\mu$ is a signed Radon measure. We aim to establish $L^p$ estimates for the gradient of $u$ that are independent of the choice of the metric $g$. This is particularly relevant when the complex structure approaches the boundary of the moduli space. To this end, we consider the metric $g' = e^{2u} g$ as a metric of bounded integral curvature. This metric satisfies a so-called quadratic area bound condition, which allows us to derive gradient estimates for $g'$ in local conformal coordinates. From these estimates, we obtain the desired estimates for the gradient of $u$.

\end{abstract}
\maketitle

\section{Introduction}

Let $f_k$ be a sequence of conformal immersions from the 2-dimensional disk $D$ into $\mathbb{R}^n$ satisfying $\int_D |A_k|^2 < 4\pi - \tau$ for some $\tau > 0$, where $A$ denotes the second fundamental form. Define the induced metric by $g_k = f_k^*(g_{\mathbb{R}^n}) = e^{2u_k}g_{euc}$. The celebrated H\'elein convergence theorem asserts that if the areas of $f_k(D)$ are bounded, then $\|u_k\|_{L^\infty(D_{\frac{1}{2}})}$ remains bounded and $f_k$ converges to a $W^{2,2}$-conformal immersion, or $u_k \rightarrow -\infty$ and $f_k$ collapses to a single point. H\'elein's theorem plays a crucial role in variational problems involving the Willmore functional. However, we can not directly get a meaningful limit in the collapsing case. When collapse occurs, one typically finds suitable translations $y_k$ and scalings $\lambda_k$ and studies the convergence of $\lambda_k(f_k - y_k)$. However, the areas potentially tend to infinity, so we need a gradient estimate to ensure that the areas remain locally bounded.

 In \cite{K-L}, the authors employed the gradient estimate $|\nabla_x G(x,y)| \leq \frac{C}{d(x,y)}$ to derive gradient bounds for solutions to equations of the form $-\Delta_g u = f$ on closed surfaces. However, the constant $C$ depends on the metric $g$,  rendering this approach unsuitable for sequences whose induced conformal classes diverge in moduli space.

In \cite{Li}, the first author utilized the conformal invariance of the Willmore energy to establish gradient estimates under the assumption that $f_k(D)$ extends to a closed immersed surface satisfying $\|A\|_{L^2} < C$. Several years later, in \cite{Li-Wei-Zhou}, the authors showed that such gradient estimates also hold under the density condition
$$
\frac{\vol(f_k(D)\cap B_r^n(y))}{\pi r^2}<C
$$
for sufficiently small $r$ and all $y \in f_k(D)$. By Simon's monotonicity inequality(\cite{Willmore}), the extension assumption implies this density condition.

In \cite{Li-Tang-Sun},the authors proved an intrinsic version of these estimates. In this paper, we apply the method in \cite{Li-Tang-Sun} to provide  gradient estimates for solutions to the equation
\begin{equation}\label{equ.u}
-\Delta_{g} u=\mu,
\end{equation}
where $\mu$ is a signed Radon measure.   

Equations of this type arise naturally and frequently. For instance, if 
$f\in L^1(\Sigma,g)$, then $fV_g$ can be regarded as a signed Radon measure. More importantly, the Green's function may be viewed as a solution to such an equation.

We say that $u$ is a solution to \eqref{equ.u} if $u\in L^1(\Sigma,g)$ and 
\begin{equation}\label{def.equ}
-\int_\Sigma u\Delta_g \varphi dV_g=\int_\Sigma \varphi d\mu
\end{equation}
holds for any $\varphi\in \mathcal{D}(\Sigma)$, where $\mathcal{D}(\Sigma)$ denotes the set of compactly supported smooth functions on $\Sigma$. It follows easily from Brezis-Merle's result (see Theorem \ref{B-M} in section 2) and Weyl's lemma (
cf.  \cite[Theorem 2.3.1]{Morrey}) that $u\in W^{1,p}_{loc}$
for any $p\in[1,2)$.
Since $\int_\Sigma u\Delta_g\varphi dV_g$ is independent of $g$, actually the definition of equation \eqref{equ.u} does not depend on $g$.
In fact, \eqref{def.equ} is equivalent to
$$
-2\sqrt{-1}\int_\Sigma u\partial\bar\partial \varphi=\int_\Sigma\varphi d\mu,
$$
i.e. the current $-2\sqrt{-1}\partial\bar\partial u$ is represented by the Radon measure $\mu$. From this perspective, it is natural to expect an estimate of $u$ that is independent of  $g$.

To state our estimates, we introduce the following notation. 
We set $\mathfrak{C}$ to be the set of conformal embeddings from $D$ into $\Sigma$ and
define
$$
E_{p}(u)=\sup\{\|\nabla u\circ\varphi\|_{L^p(D_\frac{1}{2})}:\varphi\in\mathfrak{C}\},
$$
where $\nabla$ denotes the standard gradient operator in the Euclidean space.

When the complex structure approaches the boundary of the moduli space,  it becomes necessary  to consider estimates on collars.
For  $0\leq a<\frac{1}{\pi}\log 2$,  we let $\mathfrak{C}_a$ be the set of conformal embeddings from $D\setminus \overline{D_a}$ into $\Sigma$ and define
$$
E_{a,p}(u)=\sup\{\|\nabla u\circ\varphi\|_{L^p(D_\frac{1}{2}\setminus \overline{D_{2a}})}:\varphi\in\mathfrak{C}_a\}.
$$
By convention, we set  $D_0=\{0\}$  and  define
$E_{a,p}(u)=0$ when $\mathfrak{C}_a=\emptyset$.

Our main result is the following:
\begin{thm}\label{main-gradient}
Let $\Sigma$ be a closed Riemann surface of genus $\genus$ and $u$ be a solution to  \eqref{equ.u}. Then for any $p\in[1,2)$
and $0\leq a<\frac{1}{\pi}\log 2$, 
$$
E_p(u)\leq C(p,\mathfrak{p})|\mu|(\Sigma),\s E_{a,p}\leq C(a,p,\genus)|\mu|(\Sigma).
$$
\end{thm}

\begin{rem}
Such a result does not hold locally. For example, on a disk, for the harmonic functions $u_k=kx^1$ we have
$\|\nabla u_k\|\rightarrow+\infty$ uniformly. 
\end{rem}

\begin{rem}
Applying the Collar Theorem, from Theorem \ref{main-gradient} we obtain a new proof of Theorem 0.1 in \cite{L-R}.
\end{rem}

The proof of Theorem \ref{main-gradient} is geometric in nature. Choose a smooth  metric $g$ with $\int_\Sigma|K_g|dV_{g}<C$ and set
$g'=e^{2u}g$. If we assume $|\mu|(\Sigma)=1$, then $g'$ is a metric with bounded integral curvature in the sense of Alexandrov (we refer to Section 2 for the definition).
By a result in \cite{Chen-Li}, $g'$ has quadratic area bound,  which yields the desired estimate.

As a corollary, we obtain an estimate on a constant curvature surface involving the metric.
\begin{cor}\label{gradient-metric}
Let $(\Sigma,g)$ be  a closed surface with $K_g=0$ or $-1$, and $u$ be a solution of \eqref{equ.u}.
Then there exists  a constant $C$ depending only  on the genus, such that
\begin{equation}\label{eq.gradient}
\int_{B_r^g(x_0)}|\nabla_g u|dV_g<Cr|\mu|(\Sigma).
\end{equation}
\end{cor}

For $p\in(1,2)$, our method does not provide a uniform estimate---the constant
$C$ depends on certain geometric data: when $K_g=0$,
\begin{equation}\label{gradient.torus}
r^{\frac{p-2}{p}}\|\nabla_g u\|_{L^p(B_r^g(x_0))}<C(p)\left(\frac{\pi r^2}{\vol(B_r^g(x_0))}\right)^{\frac{p-1}{p}}|\mu|(\Sigma);
\end{equation}
and when $K_g=-1$, 
\begin{equation}\label{gradient.global}
\|\nabla_g u\|_{L^p(\Sigma,g)}<C(p) \frac{|\chi(\Sigma)|^\frac{1}{p}}{{\rm Inj}(\Sigma,g)}|\mu|(\Sigma),
\end{equation}
where ${\rm Inj}$ is the injectivity radius.

\section{Preliminaries}
In this section, we list some geometric and analytic results that will be used later.
\subsection{Brezis-Merle's estimates}
The following theorem, established by Brezis and Merle ~\cite{Brezis-Merle},  plays a crucial role in this paper.
\begin{pro}[Brezis-Merle]\label{B-M}  
Given a signed Radon measure $\mu$ supported in $D\subset \R^2$ with $0<|\mu|(\R^2)<+\infty$, let 
$$
I_\mu(x)=-\frac{1}{2\pi}\int_{\R^2}\log|x-y|d\mu(y).
$$
Then $I_\mu \in W^{1,q}_{loc}(\R^2)$ for any $q\in[1,2)$ and weakly solves the 
 equation: 
\begin{equation}\label{equation.I}
-\Delta I_\mu=\mu.
\end{equation}
Moreover, we have
\begin{equation}\label{nabla.q}
r^{q-2}\int_{D_r(x)}|\nabla I_\mu|^qdx\leq C(q)|\mu|(\R^2)^q,\s \forall x, r,
\end{equation}
and
\begin{equation}\label{whynot}
\int_{D_{R}} e^{\frac{(4\pi-\epsilon)|I_\mu|}{|\mu|(\R^2)|}}dx\leq CR^\frac{\epsilon}{2\pi},\s \forall R>0\s and\s \epsilon\in(0,4\pi).
\end{equation}
\end{pro}

\begin{cor}\label{MTI} 
Let $\mu$ be a signed Radon measure on $D$ with $|\mu|(D)<\tau$. Suppose that  $u$ 
solving \eqref{equ.u} weakly on $D$ and $\|u\|_{L^1(D)}<\gamma$.    
Then for any $p<\frac{4\pi}{\tau}$ there exists  $\beta=\beta(\tau,p,\gamma)$ such that
$$
\int_{D_{\frac{1}{2}}}e^{p|u|}dx\leq \beta.
$$
\end{cor}
The proof of this corollary  can be found in \cite{Chen-Li}.

\subsection{BIC metrics} We begin with the following definition.

\begin{defi}
Let $(\Sigma,g_0)$ be a smooth Riemann surface without boundary.  
Let $\mathcal{M}(\Sigma,g_0)$ denote the set of measurable tensor $g=e^{2u}g_0$ with  $u\in L^1_{loc}(\Sigma)$ such that
 there exists a  signed  Radon measure $\mu$ satisfying 
\begin{equation*}\label{def}
\int_\Sigma\varphi \,d\mu =\int_\Sigma \left(\varphi \,K_{g_0} - u \Delta_{g_0}\varphi\right)dV_{g_0},\s \forall \varphi\in \mathcal{D}(\Sigma).
\end{equation*}
We call  $\mu$ the   Gauss curvature measure of $g$, and denote it by
$\K_g$. 
\end{defi}

Given $g\in\mathcal{M}(\Sigma,g_0)$, we  define 
\begin{equation}\label{distance-definition}
d_{g,\Sigma}(x,y)=\inf\left\{\int_\gamma e^u ds_{g_0}: \mbox{$\gamma$ is a piecewise smooth curve from $x$ to $y$ in $\Sigma$}\right\}
\end{equation}
where $\gamma$ is parametrized by its $g_0$-arclength. The function $d_{g,\Sigma}$ is well-defined distance function whenever  $d_g(x,y)$ is finite for any $x$ and $y$. In this case, we say $g$ is a metric  of {\it bounded integral curvature} in the sense of Alexandrov. For other equivalent definitions,  see ~\cite{Troyanov1}.

If
$|\K_g|(\{x\})<2\pi$
for any $x\in\Sigma$, $d_g$ is finite  (see ~\cite{Chen-Li,R, R1, R2} for more details).

The following theorem proved in \cite{Chen-Li} The following estimate is essential for the results developed in the next section.
\begin{thm}\label{vol}
Let $(\Sigma,g_0)$ be a closed surface and $g=e^{2u}g_0\in \mathcal{M}(\Sigma,g_0)$ with $|\K_g|(\Sigma)<+\infty$. Assume $d_{g,\Sigma}$ is finite in $\Sigma\times\Sigma$. Then 
\begin{equation}\label{vol.com}
\frac{\vol(B_R^g(x))}{\pi R^2}\leq 1+\frac{1}{2\pi}\K_g^-(\Sigma).
\end{equation}
\end{thm}

\section{The gradient estimate}

We begin by proving a gradient estimate on the two-dimensional disk. As mentioned in the introduction, the examples $u_k = k x^1$ show that no gradient estimates can hold when $k$ is sufficiently large. Therefore, such examples must be excluded.
\begin{defi}
We say $g=e^{2u}g_{euc}\in\mathcal{M}(D)$ has quadratic area bound with constant $\Lambda$, if for any
$x$ and $r$, we have
$$
\frac{\vol(B_r(x,d_g))}{\pi r^2}<\Lambda.
$$
\end{defi}

\begin{rem}
Suppose $\Omega_1\subset\subset\Omega_2$ and $g\in\mathcal{M}(\Omega_2)$. Then we also have $g\in\mathcal{M}(\Omega_1)$. Let $d_1$ and $d_2$ denote the distances induced  by $g$ on $\Omega_1$ and $\Omega_2$, respectively. 
Note that $d_{2}(x,y)$ may be strictly less than $d_{1}(x,y)$. However, for any $x,y\in\Omega_1$ with $d_1(x,y)< r$, we have
$$
d_2(x,y)\leq d_1(x,y)< r,
$$
which implies
$y \in B_r(x, d_2)$
, and thus \[
B_r(x, d_1) \subset B_r(x, d_2).
\]
Therefore, Theorem~\ref{vol} implies that any topological disk domain of $\Sigma$ has quadratic area bound.
\end{rem}

First, we show Theorem  1.3 in \cite{Li-Tang-Sun} still holds if we only assume that $g\in\mathcal{M}(D)$.

\begin{lem}
Let $g=e^{2u}g_{euc}\in\mathcal{M}(D)$ has quadratic area bound with constant $\Lambda$. Then for any $p$, there exist constants $C$ and $\epsilon_0$, such that if $|\K_g(D)|<\epsilon_0$, then
$$
\|\nabla u\|_{L^p(D_\frac{1}{2})}<C.
$$
\end{lem}

\proof Assume the result does not hold,  then there exists a sequence  $g_k=e^{2u_k}g_{euc}\in\mathcal{M}(D)$ such that $\|\K_{g_k}\|(D)\rightarrow 0$, while $\|\nabla u_k\|_{L^p(D_\frac{1}{2})}\rightarrow+\infty$.

Let $v_k=I_{\K_{g_k}}$ be as in Proposition \ref{B-M}. By Proposition \ref{B-M},
$$
\sup_{D_r(x)\subset D}r^{p-2}\int_{D_r(x)}|\nabla v_k|^p\rightarrow 0.
$$
Set
$w_k=u_k-v_k$, which is harmonic on $D$. Then $\|\nabla w_k\|_{C^0(D_\frac{1}{2})}\rightarrow+\infty$.

Define
$$
\rho_{k}=\sup_{x \in D_{\frac{2}{3}}}(\frac{2}{3}-|x|)|\nabla w_k|(x).
$$
There exists $x_k \in D_{\frac{2}{3}}$, such that $\rho_k=(\frac{2}{3}-|x_k|)  |\nabla w_k|(x_k)$.  Set
$r_k=1/|\nabla w_k(x_k)|$.

Since
$$
\rho_k\geq(\frac{2}{3}-\frac{1}{2})\|\nabla w_k\|_{C^0(D_\frac{1}{2})}\rightarrow+\infty,
$$
we have
$$
\frac{r_k}{\frac{2}{3}-|x_k|}=\frac{1}{\rho_k}\rightarrow 0,
$$
which implies that $D_{Rr_k}(x_k)\subset D_\frac{2}{3}$ for any fixed  $R$, when $k$
is sufficiently large.

Define $w_k'(x)=w_k(x_k+r_kx)-w_k(x_k)$.  We claim that for large $k$
\[
\sup_{x\in D_R}|{\nabla w_k}|(x) \leq 2.
\]
Indeed, for large $k$ and any $x\in D_R$, $x_k+ r_k x \in D_{\frac{2}{3}}$,  
$$
(\frac{2}{3}-|x_k|) |\nabla w_k|(x_k)\geq (\frac{2}{3}-|x_k+r_k x|)|\nabla w_k|(x_k+ r_k x),
$$
which gives
\begin{eqnarray*}
|{\nabla w_k'}|(x) &=& r_k|\nabla w_k(x_k+r_kx)| =\frac{|\nabla w_k(x_k+r_kx)|}{|\nabla w_k(x_k)|}
\leq \frac{\frac{2}{3}-|x_k|}{\frac{2}{3}-|x_k+r_kx|}
\leq \frac{\frac{2}{3}-|x_k|}{\frac{2}{3}-|x_k|-r_k R}\leq 2,
\end{eqnarray*}
for $k$ large enough.

By the standard elliptic estimate, we may pass to some subsequence and assume that $w_k'$ smoothly converges to a harmonic function $w$ on every bounded domain of $\R^2$,
such that $w(0)=0$, $|\nabla w|(0)=1$, $|\nabla w|\leq2$. 
By
 Liouville's theorem for harmonic functions,$\nabla w$ is constant, therefore $w=ax^1+bx^2$ with $a^2+b^2=1$. By rotation, we may assume $w=x^1$.

Next, we define
\[
v_k'(x) = v_k(r_kx+x_k)-c_k,\s u_k'(x)=v_k'(x)+w_k'(x),
\]
where $c_k$ is chosen  such that $\int_Dv_k'=0$. Then for any $x\in D_{R}$, we have
$$
\int_{D(x)}|\nabla v_k'|^p=r^{p-2}\int_{D_r(x_k+r_kx)}|\nabla v_k|^q\rightarrow 0.
$$
By the Poincar\'e inequality  (cf. \cite[Theorem 5.4.3]{Attouch-Buttazzo-Michaille}) and Proposition \ref{B-M}, we may assume $v_k$ weakly converges to 0 in $W^{1,q}_{loc}(\R^2)$. By Corollary \ref{MTI} and the Lebesgue-Vitali Convergence Theorem(cf. \cite[Theorem 4.5.4]{Bogachev}), for any $p\in[1,+\infty)$ and $R>0$,
\begin{equation}\label{area.dis}
\int_{D_R}|e^{u_k'}-e^{x^1}|^p\rightarrow 0,\s \int_{D_R}|e^{2u_k'}-e^{2x^1}|\rightarrow 0.
\end{equation}

Next, we compute the area bound of the metric $g_0:=e^{2x^1}g_{euc}$.
Let $T(\theta)$ be the constant such that 
$$
Length((\cos\theta,\sin\theta)t|_{t\in[0,T(\theta)]},g_0)=R,\s i.e.\s
\int_0^{T(\theta)}e^{r\cos\theta}d\theta=R.
$$
It is easy to verify that
$$
T(\theta)=\frac{\log (1+R\cos\theta)}{\cos\theta},\s i.e.\s
e^{T(\theta)\cos\theta}=R\cos\theta+1
$$
Let $a>0$ be sufficiently small and define
$$
\Omega(R)=\{(r,\theta):r\in(0,T(\theta)),\s \theta\in(-\frac{\pi}{2}+a,\frac{\pi}{2}-a)\}.
$$
Then
\begin{eqnarray*}
\vol(\Omega(R),g_0)&=&
\int_{-\frac{\pi}{2}+a}^{\frac{\pi}{2}-a}\int_0^{T(\theta)}e^{2r\cos\theta}rdrd\theta\\
&=&
\int_{-\frac{\pi}{2}+a}^{\frac{\pi}{2}-a}\left.\left(
\frac{e^{2r\cos\theta}r}{2\cos\theta}-\frac{e^{2r\cos\theta}}{4\cos^2\theta}\right)\right|_0^{T(\theta)}d\theta\\
&=&
\int_{-\frac{\pi}{2}+a}^{\frac{\pi}{2}-a}\left(\frac{(R\cos\theta+1)^2T(\theta)}{2\cos\theta}-\frac{(R\cos\theta+1)^2-1}{4\cos^2\theta}\right)d\theta\\
&=&R^2\int_{-\frac{\pi}{2}+a}^{\frac{\pi}{2}-a}\left(\frac{1}{2}\log(1+R\cos\theta)
-\frac{1}{4}\right)d\theta
+\int_{-\frac{\pi}{2}+a}^{\frac{\pi}{2}-a}RT(\theta)d\theta\\
&&+\int_{-\frac{\pi}{2}+a}^{\frac{\pi}{2}-a}\frac{T(\theta)-R}{2\cos\theta}d\theta\\
&=&R^2\int_{-\frac{\pi}{2}+a}^{\frac{\pi}{2}-a}\left(\frac{1}{2}\log(1+R\cos\theta)
-\frac{1}{4}\right)d\theta
+O(R\log R).
\end{eqnarray*}
Therefore,
\begin{equation*}
\lim_{R\rightarrow+\infty}\frac{\vol(\Omega(R),g_0)}{\pi R^2}=+\infty.
\end{equation*}

Since $\Omega(R)\subset B_R(0,d_{g_0})$, it follows that
\begin{equation*}
\lim_{R\rightarrow+\infty}\frac{\vol(B_R(0,d_{g_0},g_0)}{\pi R^2}=+\infty.
\end{equation*}

Let $T_1=\max_{|\theta|\leq\frac{\pi}{2}-a} T(\theta)$.  
Fix $p>2$.
By the first part of \eqref{area.dis} and H\"older's inequality,
$$
\int_{0}^{2\pi}\int_{0}^{T_1}|e^{u_k'(r,\theta)}-e^{w(r,\theta)}|drd\theta=
\int_{0}^{2\pi}\int_{0}^{T_1}|e^{u_k'(r,\theta)}-e^{w(r,\theta)}|r^\frac{1}{p}r^{-\frac{1}{p}}drd\theta\leq
C\|e^{u_k'}-e^{w}\|_{L^p(D_{T_1})}\rightarrow 0.$$
Passing to a subsequence, we have $\int_{0}^{T_1}|e^{u_k'(r,\theta)}-e^{w(r,\theta)}|dr$ converges to 0 for a.e. $\theta\in S^1$. For small $\epsilon>0$, by Egorov's theorem there exists $A\subset S^1$ with Lebesgue measure less than $\epsilon$ such that  
$
\int_0^{T_1}|e^{u_k'}-e^{w}|dr\rightarrow 0$
uniformly on $S^1\setminus A$. Set $g_k'=e^{2u_k'}g_{euc}$ and
$$
\Omega(R,A)=\Omega(R)\setminus \{(r,\theta):\theta\in A\}.
$$
By the Trace Embedding Theorem, For large $k$,
$$
Length((\cos\theta,\sin\theta)t|_{t\in[0,T(\theta)]},g_k')<R+1,\s \forall \theta\in S^1\setminus A,
$$
so
$$
\Omega(R,A)\subset B_{R+1}^{g_k'}(0).
$$
By the second part of \eqref{area.dis}, for small $\epsilon>0$ and large $k$, 
$$
\vol(\Omega(R,A),g_k')\geq
\frac{1}{2}\vol(\Omega(R),g_0).
$$
Therefore, for suitable  $R>0$,  
$$
\frac{\vol(B_{(R+1)r_k}^{g_k}(x_k),{g_k})}{\pi((R+1)r_k)^2}=\frac{\vol(B_{R+1}^{g_k'}(0),g_k')}{\pi(R+1)^2}>\Lambda,
$$
which contradicts the quadratic area bound condition.
\endproof

Next, we remove the assumption that $|\K_g|(D)<\epsilon_0$.

\begin{thm}\label{general.estimate}
Let $g=e^{2u}g_{euc}\in\mathcal{M}(D)$. Assume  $g$ has quadratic area bound with constant $\Lambda$ and $|\K_g|(D)<\Lambda'$. Then for any $p\in [1,2)$, there exists a constant $C=C(\Lambda,\Lambda',p)$  such that 
$$
\|\nabla u\|_{L^p(D_\frac{1}{2})}<C.
$$
\end{thm}

\proof
Suppose the conclusion fails, then there exists a sequence $g_k=e^{2u_k}g_{euc}$ with $|\K_{g_k}|(D)<\Lambda'$
and quadratic area bound with constant $\Lambda$ such that
$$
\|\nabla u_k\|_{L^p(D_\frac{1}{2})}\rightarrow+\infty.
$$
Let $v_k=I_{\K_{g_k}}$ and $w_k=u_k-v_k$.
By Proposition \ref{B-M},
$$
r^{p-2}\int_{D_r(x)}|\nabla v_k|^p<C,\s \forall D_r(x)\subset D.
$$
Thus,
$$
\sup_{D_\frac{1}{2}}|\nabla w_k|\rightarrow+\infty.
$$

Passing to a subsequence if necessary, we may assume that $|\K_{g_k}|$ converges weakly to a Radon measure $\mu$. Set
$$
\mathcal{S}=\{x\in D: \mu(\{x\})>\frac{\epsilon_0}{2}\}.
$$
Clearly, $\mathcal{S}$ is a finite set. For any $x\notin \mathcal{S}$, we can find $r>0$, such that $|\K_{g_k}|(D_r(x))<\epsilon_0$ for sufficiently large $k$.
Then  $\|\nabla u_k\|_{L^p(D_{r/2}(x))}<C(r)$, which implies that $\|\nabla w_k\|_{L^p(D_{r/2}(x))}<C(r)$. By the standard elliptic estimate, after passing to a subsequence we may assume that there exist constants $c_k$(for instance, fix a point $p_0\notin\mathcal{S}$ and let $c_k=w_k(p_0)$) and a harmonic function $w$, such that $w_k-c_k\to w$ smoothly on $\Omega$ for any $\Omega\subset\subset D\setminus\mathcal{S}$. 

Choose $r_0>0$ such that for any distinct $p,p'\in\mathcal{S}$,  $D_{r_0}(p)\cap D_{r_0}(p')=\emptyset$. For any $x\in D_{r_0/4}(p)$, $\partial D_{r_0/2}(x)\subset D_{r_0}(p)\setminus D_{r_0/4}(p)$.
Applying the mean value property of harmonic functions on
$\partial D_{r_0/2}(x)$, it turns out that $w_k-c_k$ are uniformly bounded on $D_{r_0/4}(p)$. Therefore,
$|\nabla w_k|$ are uniformly bounded on $D_{r_0/4}(p)$, hence $|\nabla w_k|$ are uniformly bounded on $D_{\frac{1}{2}}$,
which leads to a contradiction.
\endproof

{\it Proof of Theorem \ref{main-gradient}:}  
Without loss of generality, we assume $|\mu|(\Sigma)=1$.
Choose a smooth metric $g$ with $\int_\Sigma|K_g|dV_g<C(\genus)$, where $C(\genus)$ depends only on the genus $\genus$ of $\Sigma$. Set $g'=e^{2u}g\in\mathcal{M}(\Sigma,g)$. Then 
$$
\K_{g'}=\mu-K_gdV_g.
$$
Moreover, we have $|\K|_{g'}(\{x\})<2\pi$ for any $x$, so $d_{g'}$ is finite. By Theorem \ref{vol},   both $g'$ and $g$ have  
quadratic area bound with some constant $\Lambda>0$.

For any $\varphi\in\mathfrak{C}$, $\varphi$ defines a local isothermal coordinate chart. In this chart, we write 
\[
\varphi^*(g) = e^{2v} g_{\mathrm{euc}}, \quad \varphi^*(g') = e^{2(u + v)} g_{\mathrm{euc}}.
\]
By Theorem \ref{general.estimate},  
$$
\|\nabla v\|_{L^p(D_\frac{1}{2})}<C,\s \|\nabla (u+v)\|_{L^p(D_\frac{1}{2})}<C
$$
which implies the first part of Theorem \ref{main-gradient}. \\

Next, we consider $\varphi\in\mathfrak{C}_a$ and define $v$ as above. 
For simplicity, let $a=2^{-m}$ with $m>2$ being a positive integer. Note that for any $x\in D_{2^{2-j}}\setminus \overline{D_{2^{1-j}}}$(where $2<j\leq m$), we have $D_{2^{-j}}(x)\subset D\setminus\overline{D_a}$.
By a rescaling argument, 
$$
\int_{D_{2^{-j}}(x)}|\nabla u|^p\leq C(2^{-j})^{2-p}.
$$
By a covering argument, 
$$
\int_{D_{2^{-i}}\setminus \overline{D_{2^{-i-1}}}}|\nabla u|^p\leq C(2^{-i})^{2-p},
$$
so
$$
\int_{D_{2^{-i}\setminus D_{2^{-m}}}}|\nabla u|^pdx\leq C\sum_{i=1}^m (2^{-i})^{2-p}.
$$

For the case $a=0$, letting $m\to\infty$ we get
$$
\int_{D_{2^{-1}}\setminus\{0\}}|\nabla u|^p\leq  C\sum_{i=1}^\infty (2^{-i})^{2-p}.
$$
\endproof

\section{Estimates on a surface of constant curvature}
In this section, we assume that $K_g=0$ or $-1$ on $\Sigma$ and that $u$ is a solution of \eqref{equ.u}.  The main aim  is to prove Corollary \ref{gradient-metric}, as well as \eqref{gradient.torus} and \eqref{gradient.global}. For simplicity, we assume $|\mu|(\Sigma)=1$ in this section.

\subsection{The flat case}
We first consider the case $K_g=0$.  
 Note that \eqref{eq.gradient} and \eqref{gradient.torus} hold for $g$ if and only if they hold for $\lambda g$ for any $\lambda>0$.
 Without loss of generality, we assume $(\Sigma,h)$ is induced by the lattice
$\{1,a+b\sqrt{-1}\}$ in $\C$, where $-\frac{1}{2}<a\leq\frac{1}{2}$, $b>0$, $a^2+b^2\geq 1$, and
$a\geq 0$ whenever $a^2+b^2=1$.  Set 
$$
\rho=\sqrt{a^2+b^2},\s \theta=\arccos \frac{a}{\rho},\s
v=(\rho,0).\s   w=(\cos\theta,\sin\theta).
$$
Then $(\Sigma,h)$ is also induced by the lattice $\{v,w\}$,
where $\rho\geq1$ and $\frac{\pi}3\leq\theta<\frac{2\pi}3$. 

Let $\Pi: \mathbb{C} \to \Sigma$ be the covering map. 
Let $u'$ denote the lift of $u$, i.e. $u'=u\circ \Pi$.
 For any $x\in\C$, and $r<\frac{\sqrt{3}}{4}$, $D_r(x)$ can be viewed as an isothermal coordinate chart of $\Sigma$
around $\Pi(x)$. By Theorem \ref{main-gradient}, 
$$
 \int_{D_r(x)}|\nabla u'|^pdx\leq Cr^{2-p}.
$$

By a covering argument, we have
$$
\int_{[i,i+1]\times[-\frac{\sin\theta}{2},\frac{\sin\theta}{2}]} |\nabla u'|^p<C.
$$

Without loss of generality, we assume $0$ is a lift of $x_0$. 
When $r<3$,   note that $r/8<\frac{\sqrt{3}}{4}$. 
Covering 
$D_r(0)$ with finitely many $D_\frac{r}{8}(x)$
(having a universal upper bound on the number of disks) gives
 $$
\int_{B_r^g(x_0)}|\nabla_g u|^pdV_g\leq  \int_{D_r(0)}|\nabla u'|^pdx\leq Cr^{2-p}.
 $$
When $p=1$, this yields \eqref{eq.gradient}. When $p\in(1,2)$, since
$0<C\leq \frac{\vol(B_r^g(x_0))}{\pi r^2}\leq 1$, we obtain \eqref{gradient.torus}.

When $r>3$,   let $m$ be the smallest integer such that 
$r\leq m$. It is easy to check that if $x\in \Pi( [-\frac{m}{2},\frac{m}{2}]\times[-\frac{\sin\theta}{2},\frac{\sin\theta}{2}])$, then
$$
d_g(x,x_0)\leq \sqrt{(\frac{m}{2})^2+(\frac{\sin\theta}{2})^2}\leq m-1\leq r,
$$
so
\begin{equation}\label{torus.area}
\Pi( [-\frac{m}{2},\frac{m}{2}]\times[-\frac{\sin\theta}{2},\frac{\sin\theta}{2}])\subset B_r(x_0)\subset \pi([-m,m]\times[-\frac{\sin\theta}{2},\frac{\sin\theta}{2}] ),
\end{equation}
Then
$$
\int_{B_r(p)}|\nabla_gu|^pdV_g\leq \int_{[-m,m]\times[-\frac{\sin\theta}{2},\frac{\sin\theta}{2}]}|\nabla u'|^p=
\sum_{i=-m}^{m-1}\int_{[-m,m]\times[-\frac{\sin\theta}{2},\frac{\sin\theta}{2}]}|\nabla u'|^p\leq Cm.
$$
When $p=1$, using $m-1< r$, we get \eqref{eq.gradient}. When $p\in(1,2)$, by \eqref{torus.area},
$$
\vol(B_r^g(x_0))\leq Cm<2Cr,
$$
so
$$
\frac{r^p}{\vol^{p-1}(B_r^g(x_0))}\geq Cr^p/r^{p-1}=Cr,
$$
which yields \eqref{gradient.torus}.

\subsection{On a hyperbolic surface}
Now, we assume  $K_g=-1$. 

Let  $r_x$ denote the injectivity radius of $g$ at $x$. Since 
$$
\vol(B_{r_x}^g(x))=2\pi(\cosh r_x-1)\leq \vol(\Sigma)=2\pi|\chi(\Sigma)|,
$$
it follows that $r_x$ is bounded above by a constant depending only on the genus of $\Sigma$. 

It is well-known that for any $r<r_x$
$(B^g_{r}(x),g)$ is biholomorphic to $(D,g_r)$, where
$$
g_r=\left(\frac{2\sinh(r/2)}{1-\sinh^2(r/2)|x|^2}\right)^2g_{euc}.
$$
By applying Theorem \ref{main-gradient}, one readily verifies that  when $r\leq r_x/2$,  
$$
\int_{B^g_r(x)}|\nabla_g u|^pdV_g\leq Cr^{2-p}.
$$

Let $a>0$ be fixed and $\Sigma_{a}$ be the set of points where $r_x\geq a$. By Vitali's 5-times covering theorem, there exist pairwise disjoint closed balls $\{\overline{B^g_{a/10}(x_i)}\}$ with $x_i\in\Sigma_{a}$
such that
$\{\overline{B_{a/2}(x_i)}\}$ cover $\Sigma_{a}$. Then 
\[
\# \{ x_i \} \le \frac{ \vol(\Sigma) }{ V_{a/10} },
\]
 where $V_{a/10}$ is the area of a disk of radius $a/10$ in the hyperbolic plane. Therefore
$$
\int_{\Sigma_{a}}|\nabla_gu|^p<C\frac{|\chi(\Sigma)|}{a^2}\left(\frac{a}{2}\right)^{2-p}<C\frac{|\chi(\Sigma)|}{a^p},
$$
which gives \eqref{gradient.global}.

Next, we consider  the collar regions, which are components of $\Sigma\setminus\Sigma_{\arcsinh 1}$.  For more details on collars, we refer the reader to \cite{Buser}.

Let $\gamma$ be a geodesic loop of length $\ell(\gamma)<2\arcsinh 1$,  $\gamma$ defines a collar  $U\subset\Sigma$,  isometric to
$ S^1\times(-T,T)$ with the metric 
$$
\left(\frac{\lambda}{\cos\lambda t}\right)^2(dt^2+d\theta^2).
$$
See the appendix for the definitions of $T$ and $\lambda$.

We claim that for small $\lambda>0$ and any $t_1,t_2$ satisfying $|t_1|<T-1$,$|t_2|<T-1$, $|t_2-t_1|<2$,  there exists a universal constant $C_0>1$ such that
\begin{equation}\label{ratio}
C_0^{-1}<\frac{\cos\lambda t_2}{\cos\lambda t_1}<C_0    
\end{equation}
In fact,\eqref{ratio} holds with $C_0=e^2$.
Indeed, Lagrange's mean value theorem (together with the fact $\lambda(T-1)<\frac{\pi}2-\lambda$) gives that
$$
|\log\cos\lambda t_2-\log\cos\lambda t_1|\leq\lambda|t_2-t_1|\sup_{|t|<\frac{\pi}2-\lambda}|\tan t|<\frac{2\lambda}{\tan\lambda}<2
$$
From which the claim follows.

Let $k,m$ be integers.We claim that
$$
\int_{[k,m]\times S^1}|\nabla_g u|dV_g\leq Cd_{k,m}
$$
for any $-T+2<k<m<T-2$, where 
we define 
\begin{equation}\label{dist}
    d_{t,t'}=d_g(\{t\}\times S^1,\{t'\}\times S^1).
\end{equation}
We first prove the case $k\geq 0$ . 
Let $Q_i=[i,i+1]\times S^1$. By Theorem \ref{main-gradient},  
\begin{equation}\label{tube}
    \int_{Q_i}|\nabla u|dtd\theta<C.
\end{equation}
 
Since $t\mapsto\frac{\lambda}{\cos\lambda t}$
is increasing on $[k,m]$, \eqref{tube} implies that
\begin{eqnarray*}
\int_{[k,m]\times S^1}|\nabla_g u|dV_g\leq C\sum_{i=k}^{m-1}\frac{\lambda}{\cos \lambda(i+1)}\leq C
\int_{k+1}^{m+1}\frac{\lambda}{\cos \lambda s}ds
.
\end{eqnarray*}
Note that, by \eqref{ratio}  and the monotonicity of $\cos t$ on $[0,\frac{\pi}{2}]$,
\begin{equation}\label{length.i.i+1}
\int_{i-1}^{i}\frac{\lambda ds}{\cos \lambda s}\leq \int_{i}^{i+1}\frac{\lambda ds}{\cos\lambda s}\leq C\int_{i-1}^{i}\frac{\lambda ds}{\cos \lambda s}.
\end{equation}
Hence
$$
\int_{k+1}^{m+1}\frac{\lambda ds}{\cos\lambda s}\leq C\int_{k+1}^{m-1}\frac{\lambda ds}{\cos\lambda s}\leq Cd_{k,m}.
$$
Thus
$$
\int_{[k,m]\times S^1}|\nabla_g u|dV_g\leq Cd_{k,m}.
$$

When $k< 0$, we have
$$
\int_{[k,m]\times S^1}|\nabla_g u|dV_g=\int_{[k,0]\times S^1}|\nabla_g u|dV_g+\int_{[0,m]\times S^1}|\nabla_g u|dV_g\leq C(d_{k,0}+d_{0,m})=Cd_{k,m}.
$$
This proves our claim.\\

Next, we estimate $\int_{(B_r^g(x_0)\setminus\Sigma_a)\cap U}|\nabla u|$.
Choose $a>0$ sufficiently small
and assume $(B^g_r(x_0)\setminus\Sigma_a)\cap U\neq\emptyset$, then   $\frac{1}{2}\ell(\gamma)=\min_{x\in U} r_x<a$. We can choose $a$ such that (see the calculations in the appendix)
$$
\inf_{t\in(-T,-T+10)\cup (T-10,T)} r_{(t,\theta)}>a,\s 
and\s d_{T-10,T}>2a.
$$

We have two cases:

\emph{Case 1.} $B_r^g(x_0)\subset[k-1,k+1]\times S^1\subset(-T+5,T-5)\times S^1$. Without loss of generality, we assume $x_0=(t_0,0)$.

{It suffices to consider the case where $r>r_{x_0}/2$. Then there exists a unit-speed geodesic $\gamma:[0,2r_{x_0}]\rightarrow B_{2r}^g(x_0)$ with $\gamma(0)=\gamma(2r_{x_0})$. In particular, 
$\gamma$ must pass through a point  $x_0'=(t_0',\pi)$. By \eqref{ratio}, we have
$0< \cos\lambda t\leq C\cos\lambda t_0$, which implies 
$$
\frac{\lambda}{\cos\lambda t_0}\leq C d(x_0',x_0)\leq 2Cr.
$$
Consequently, 
$$
\int_{B_r^g(x_0)}|\nabla_gu|dV_g\leq \frac{\lambda}{\cos\lambda t_0}\int_{[t_0-1,t_0+1]\times S^1}|\nabla u|dtd\theta\leq C\frac{\lambda}{\cos\lambda t_0}\leq Cr.
$$
}

\emph{Case 2.} Case 1 does not hold. Let $t_1$ and $t_2$ be the smallest and largest $t\in[-T+5,T-5]$ respectively, such that $\{t\}\times S^1\cap \overline{B_r^g(x_0)}\neq\emptyset$. Then 
$$
d_g(\{t_1\}\times S^1,\{t_2\}\times S^1)\leq 2 r.
$$
Let $m=[t_2]+1$, $k=[t_1]-1$. By \eqref{length.i.i+1}, 
$$
d_{k,m}\leq Cd_{t_1,t_2}\leq Cr.
$$
Thus
$$
\int_{B_r(x_0)\setminus\Sigma_a\cap U}|\nabla_g u|dV_g\leq Cd_{k,m}\leq Cr.
$$

We are now in a position to prove \eqref{eq.gradient}.  
When $r<a$, we  have already proved it for $x_0$ in either $\Sigma_a$ or $\Sigma \setminus \Sigma_a$.   When $r\geq a$,
since $\Sigma$ has at most $3\genus -3$ collars, we have
$$
\int_{B_r^g(x_0)}|\nabla_gu|dV_g=
\int_{B_r^g(x_0)\cap\Sigma_a}|\nabla_gu|dV_g+\int_{B_r^g(x_0)\setminus\Sigma_a}|\nabla_gu|dV_g\leq C+ Cr\leq C'r.
$$
~\endproof

\section{Appendix}
In this appendix, we list some facts about collars.

Let $\gamma$ be a geodesic loop with length $\ell(\gamma)<2\arcsinh 1$, and
$$
w= \arcsinh\frac{1}{\sinh(\frac{1}{2}\ell(\gamma))}.
$$
By \cite[Theorem 4.1.6]{Buser}, $\gamma$ defines a collar  $U\subset\Sigma$,  isometric to
$ S^1\times(-w,w)$ with the metric 
$$
g=d\rho^2+ \ell^2(\gamma)\cosh^2\rho ds^2=\left(\frac{\ell(\gamma)\cosh\rho}{2\pi}\right)^2\left(\left(\frac{2\pi d\rho}{\ell(\gamma)\cosh\rho}\right)^2+d\theta^2\right),
$$ 
where $s=\frac{\theta}{2\pi}$.

Let 
$$
(t,\theta)=\phi(\rho,\theta)=\left(\frac{4\pi\arctan e^\rho}{\ell(\gamma)},\theta\right).
$$
Then $\phi$ is a diffeomorphism from $(-w,w)\times S^1$
to 
$$
\left(\frac{4\pi\arctan e^{-w}}{\ell(\gamma)},\frac{4\pi\arctan e^{w}}{\ell(\gamma)}\right)\times S^1
$$
with
$$
\phi*(g)=\ell^2(\gamma)\cosh^2\rho(dt^2+d\theta^2)=\left(\frac{\ell(\gamma)}{2\pi\sin\frac{\ell(\gamma)t}{2\pi}}\right)^2(dt^2+d\theta^2).
$$

Set
$$
T=\frac{4\pi\arctan e^{w}}{\ell(\gamma)}-\frac{\pi^2}{\ell(\gamma)},\s \lambda=\frac{\ell(\gamma)}{2\pi}.
$$
Then the collar is also isometric to 
$$
\left((-T,T)\times S^1,\left(\frac{\lambda}{\cos\lambda t}\right)^2(dt^2+d\theta^2)\right).
$$
For $-T\leq t<t'\leq T$, recall that we have defined (see\eqref{dist}) 
$$
d_{t,t'}=d_g(\{t\}\times S^1,\{t'\}\times S^1).
$$

We compute $d_{T-t,T}$ and $r_{(T-t,\theta)}$  for fixed $t$ and
small $\ell:=\ell(\gamma)$. 
Recall the asymptotics as  $x\rightarrow+\infty$
$$
\arctan x=\frac{\pi}{2}-\frac{1}{x}+O(\frac{1}{x^2}),\s 
\arcsinh  x = \log(2x) + O(\frac{1}{ x^2}).
$$
We have
$$
w=\arcsinh\frac{1}{\ell/2+O(\ell^3)}=\arcsinh(\frac{2}{\ell}+O(\ell))=\log\frac{4}{\ell}+O(\ell),
$$
Then it follows from Lagrange's mean value theorem that
$e^{-w}=\frac{\ell}{4}+O(\ell^2)
$,
hence
$$
\s
\arctan e^w=\frac{\pi}{2}-\frac{\ell}{4}+O(\ell^2),\s 
T=\frac{\pi^2}{\ell}-\pi+O(\ell).
$$
Thus
\begin{eqnarray*}
d(t):=d_{T-t,T}&=&\int_{T-t}^T\frac{\lambda}{\cos\lambda s}ds=\int_{\lambda(T-t)}^{\lambda T}\frac{ds}{\cos s}=\log\frac{1+\sin\lambda T}{1+\sin\lambda(T-t)}-
\log\frac{\cos \lambda T}{\cos \lambda(T-t)}\\
&=&-\log\frac{\pi}{\pi+t}+O(\ell),
\end{eqnarray*}
by \cite[Theorem 4.1.6]{Buser} again
$$
\sinh r_{(T-t,\theta)}=\cosh\frac{\ell}{2}\cosh d(t)-\sinh d(t)=e^{-d(t)}+O(\ell)=
\frac{\pi}{\pi+t}+O(\ell),
$$
so
$$
 r_{(T-t,\theta)}=\arcsinh\frac{\pi}{\pi+t}+O(\ell).
$$

\section*{Acknowledgment}
The first author is partially supported by  National Key R\&D Program of China 2022YFA1005400.

\end{document}